\documentclass[10pt]{article}

\usepackage{amsmath,amsthm,enumerate,amssymb}
\usepackage{bm}
\usepackage[english]{babel}
\usepackage{graphicx}

\newcommand{\md}{\,\mathrm{d}} 

\newcommand{\Nabla}{\nabla}

\newcommand{\C}{{\mathbb{C}}}
\newcommand{\E}{{\mathbb{E}}}

\DeclareRobustCommand\openone{\leavevmode\hbox{\small1\normalsize\kern-.33em1}}

\newcommand{\be}{\begin{equation}}
\newcommand{\ee}{\end{equation}}
\newcommand{\bea}{\begin{eqnarray}}
\newcommand{\eea}{\end{eqnarray}}
\newcommand{\beas}{\begin{eqnarray*}}
\newcommand{\eeas}{\end{eqnarray*}}

\newtheorem{theorem}{Theorem}
\newtheorem{lemma}{Lemma}

\newtheorem{corollary}{Corollary}

\newtheorem{proposition}{Proposition}
\newcount\minute
\newcount\hour
\def\currenttime{%
    \minute\time
    \hour\minute
    \divide\hour60
    \the\hour:\multiply\hour60\advance\minute-\hour\the\minute}

\begin{document}

\title{Inverse moments of univariate discrete distributions\\ via the Poisson expansion}
\author{Koenraad Audenaert}

\date{\today, \currenttime}

\maketitle

\begin{abstract}
In this note we present a series expansion of inverse moments of a non-negative discrete random variate
in terms of its factorial cumulants, based on the Poisson-Charlier expansion of a discrete distribution. 
We apply the general method to the positive binomial
distribution and obtain a convergent series for its inverse moments with an error
residual that is uniformly bounded on the entire interval $0\le p\le 1$. \\[1ex]
\textbf{Keywords:}
Inverse moments, factorial cumulant, positive binomial, Poisson expansion.
\end{abstract}


\section{Introduction}
Certain problems in statistics and in other branches of science require the calculation
of the inverse moments of a distribution, also called the negative or reciprocal moments.
In statistics, noteworthy examples are life testing problems \cite{Epstein:53} and
sampling problems with samples of random length \cite{Marciniak:99}.
Especially useful are the inverse moments of the binomial distribution.
Recent applications in quantum physics, for example, include the calculation of 
running times of certain quantum computation algorithms \cite{Znidaric:06}, 
the study of random walks on $n$-dimensional cubes \cite{Garcia:01},
and exact calculations of the confidence region of a Beta
estimator of the parameter $p$ of a binomial distribution \cite{Audenaert:08}.

While inverse moments of positive binomial variates can be calculated exactly \cite{Stephan:45},
the ensuing expressions are $(N+1)$-term summations and the calculations become complicated for large $N$.
In addition, it is far from clear from these summation formulas what the asymptotic behaviour with $N$ should be.
For these reasons there is an interest in obtaining efficient series expansions.
Inverse moments of positive binomial variates
have been studied as early as 1945 by Stephan \cite{Stephan:45}, who considered
expected value and variance of negative powers and expressed them as series
expansions of inverse factorials.
In that work an even earlier reference was made to the work of Bohlmann \cite{Bohlmann:13}, 
whose approach was to expand
the function $1/x$ in a Taylor series and take expected values of each term. This approach has recently
been revived in Ref.\ \cite{Znidaric:05}.
Based on these expansions, Grab and Savage \cite{Grab:54} have calculated tables for the first
inverse moment of positive binomial and Poisson variates.
Many other expansions have been proposed, for example expansions
in Eulerian polynomials \cite{Marciniak:99} and in
factorial powers \cite{Rempala:03}.
These expansions also work for related distributions like the inverse binomial
and Poisson distributions.
Govindajarulu \cite{Govi:63} has found recurrence relations between inverse moments of positive binomial variates
and Refs.\ \cite{Jones:86,Lew:76,Wooff:85,Zacks:80,Pittenger:90} contain various bounds on inverse moments.
More general methods, valid for any distribution and involving integrals, have been proposed in
\cite{Chao:72,Cressie:81,David:56,Jones:87,Rockower:88}.
The problem with any of these series expansions considered so far is that as concerns convergence
they do not perform equally well over the complete range $0\le p\le 1$. Some of these expansions are asymptotic
series and are divergent. We illustrate this claim with graphical examples in Section \ref{sec:numeric}.

In this paper we present a new series expansion for inverse moments, based on a very simple idea.
A discrete distribution of a non-negative discrete random variable can be approximated by a series expansion
based on the Poisson distribution \cite{Barbour:87}. The coefficients of this expansion are the
factorial cumulants of the original distribution. This expansion then induces an expansion of the inverse
moments in terms of the inverse moments of the Poisson distribution. The main point we wish to make in this paper
is that this expansion of inverse moments is a very good one because of its excellent convergence
properties, converging very rapidly and consistently throughout the interval $0\le p\le 1$.

In Section \ref{sec:barbour} we give a brief overview of the Poisson expansion method,
and in Section \ref{sec:expectation} we derive the
required expectation values of the Poisson distribution in terms of its inverse moments, which can be
calculated with standard software or alternatively via the formulas presented in
Appendix \ref{app:numerical}.
Our main result is Theorem 1, Section \ref{sec:main}. The excellent convergence properties are obvious from Figures
\ref{fig:poisson10abs} and \ref{fig:poisson10}, which show the absolute and relative error of the expansion applied
to the binomial distribution.
\section{Notations}
Throughout this paper, $Q$ denotes a Poissonian random variate with mean value $\mu$. Its probability distribution function (PDF) is given by the semi-infinite sequence
$\pi_\mu=(\pi_\mu(0),\pi_\mu(1),\ldots,\pi_\mu(k),\ldots)$
with $\pi_\mu(k) = e^{-\mu} \mu^k/k!$, for $k=0,1,\ldots$

\bigskip

Let $f$ be a semi-infinite sequence
$$
f=(f(0),f(1),\ldots,f(k),\ldots).
$$
We will always set $f(k)=0$ for $k<0$.
The \textit{forward difference operator} $\Delta$ and 
the \textit{backward difference operator} $\nabla$ are defined via
\beas
\Delta f(k) &=& f(k+1)-f(k) \\
\nabla f(k) &=& f(k)-f(k-1).
\eeas
These operators can be represented by semi-infinite matrices:
\beas
\Delta\mapsto \bm{\Delta} &=& \left(
\begin{array}{cccc}
-1 & 1  & \\
   & -1 & 1 & \\
   &    &\ddots&\ddots
\end{array}
\right) \\
\nabla\mapsto \bm{\nabla} &=& -\bm{\Delta}^T = \left(
\begin{array}{cccc}
1  &   && \\
-1 & 1 &&  \\
   & -1 & 1 & \\
   && \ddots &\ddots
\end{array}
\right).
\eeas
Higher-order difference operator $\Delta^l$ and $\nabla^l$ ($l>1$) are defined as
\bea
\Delta^l f(k) &=& \sum_{j=0}^l {l\choose j} (-1)^{l-j} f(k+j) \\
\nabla^l f(k) &=& \sum_{j=0}^l {l\choose j} (-1)^j f(k-j).
\eea
These operators are represented by the $l$-th matrix powers of $\bm{\Delta}$ and $\bm{\nabla}$:
\beas
\Delta^l &\mapsto& \bm{\Delta}^l \\
\nabla^l &\mapsto& \bm{\nabla}^l = (-1)^l (\bm{\Delta}^l)^T.
\eeas

\bigskip

The \textit{factorial cumulants} $\kappa^{(j)}$ of a discrete distribution with PDF $f$
are generated by the logarithm of the expectation value of $(1+x)^k$:
\be
\log\left(\sum_{k=0}^\infty f(k)\,(1+x)^k\right) = \sum_{j=0}^\infty \kappa^{(j)} \,x^j/j!.
\ee

\bigskip

To calculate certain inverse moments of the Poisson distribution, 
we will need the Stirling numbers of the first kind
$S_j^{(k)}$ \cite{Abramowitz:72}.
They satisfy the recurrence relation
\be
S_{j+1}^{(k)} = S_j^{(k-1)} - j S_j^{(k)},\quad j\ge k\ge 1,
\label{eq:recstir}
\ee
with $S_1^{(1)} = 1$,
and are generated by
\be
\prod_{k=0}^{n-1} (x-k) = \sum_{j=1}^n S_n^{(j)}x^j.
\label{eq:genstir}
\ee
We also need the non-central Stirling numbers of the first
kind $S_{j,l}^{(k)}$ \cite{Hsu:98,Koutras:82};
for $l=0$ they coincide with the ordinary Stirling numbers. They satisfy the recurrence
\be
S_{j+1,l}^{(k)} = S_{j,l}^{(k-1)} - (j+l) S_{j,l}^{(k)},\quad j\ge k\ge 1,
\label{eq:recshift}
\ee
with $S_{1,l}^{(1)} = 1$,
and are generated by
\be
x \prod_{k=l+1}^{n-1+l} (x-k) = \sum_{j=1}^n S_{n,l}^{(j)}x^j.
\label{eq:genshift}
\ee
For $k=1$ explicit formulas exist:
\be
S_j^{(1)} =(-1)^{j-1}(j-1)!
\label{eq:stirling1}
\ee
and
\be
S_{j,l}^{(1)} =(-1)^{j-1}(j+l-1)!/l!
\label{eq:stirlingnc1}
\ee

\section{Numerical comparison of existing expansions\label{sec:numeric}}
Here we consider three existing expansions of the first inverse moment of the binomial distribution
and illustrate their convergence behaviour over the entire interval $0\le p\le 1$.
As is customary, we put $q=1-p$.

\begin{figure}
\includegraphics[width=16cm]{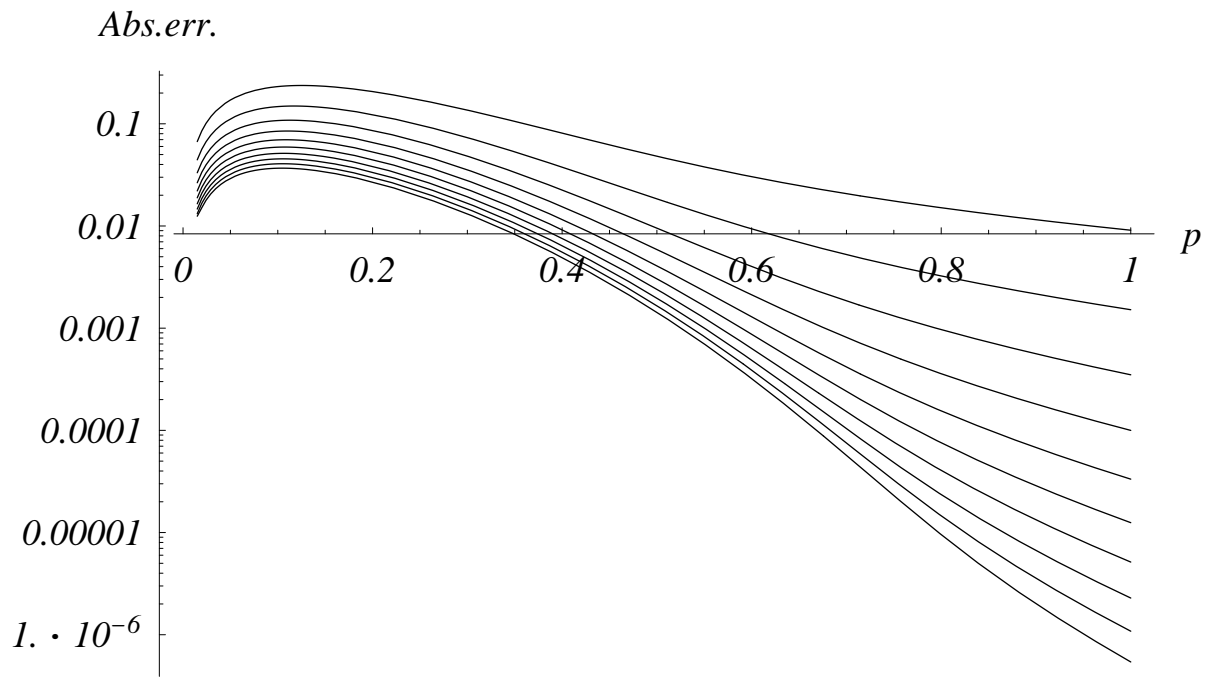}
\includegraphics[width=16cm]{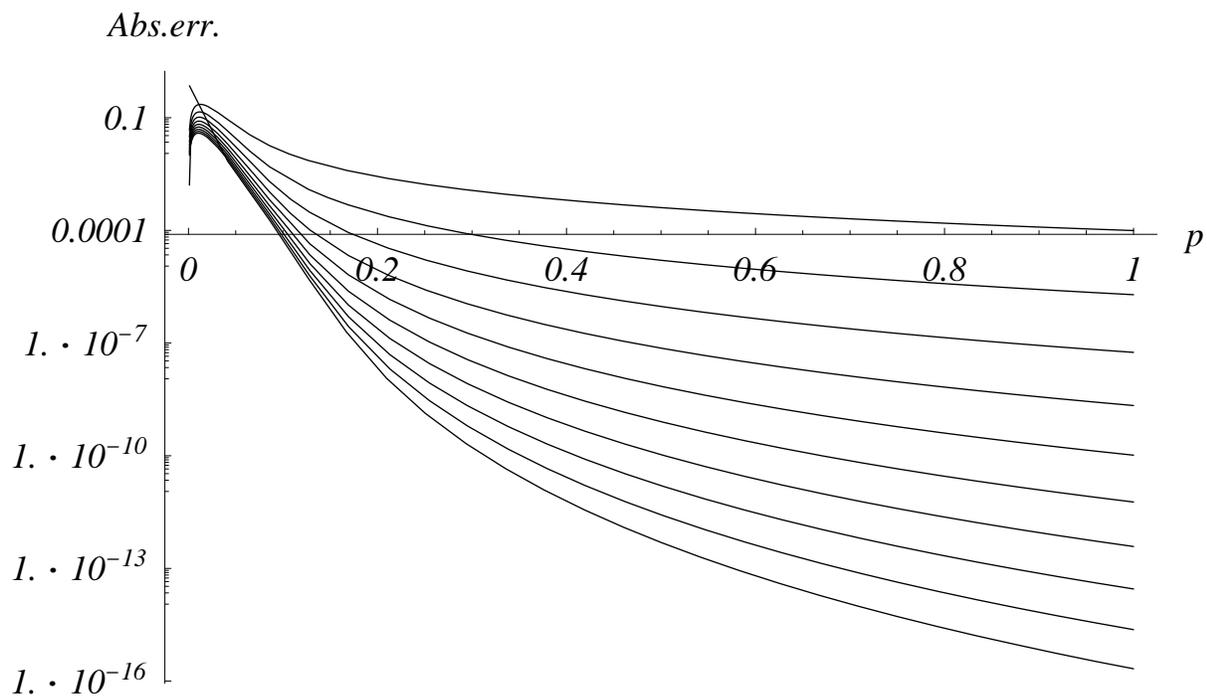}
\caption{
Absolute error as a function of $p$ of Stephan's expansion of
the first inverse moment of a positive binomial variate for $N=10$ (upper graph)
and $N=100$ (lower graph), with 1 term (upper curve), and
up to 10 terms (lowest curve).
\label{fig:stephan10}
}
\end{figure}

The oldest known expansion is Stephan's expansion \cite{Stephan:45}.
The $M$-term expansion reads
\bea
(1-q^N) \E[1/K|K>0] &=& \sum_{i=1}^M \frac{(i-1)!N! s_i}{(N+i)!p^i} \\
s_i &:=& 1-\sum_{j=0}^i {N+i\choose j} p^j q^{N+i-j}.
\eea
For very small $p$, this expansion suffers from numerical instability.
For moderately small $p$, convergence is very slow, as can be seen from Figure \ref{fig:stephan10}.

\begin{figure}
\includegraphics[width=16cm]{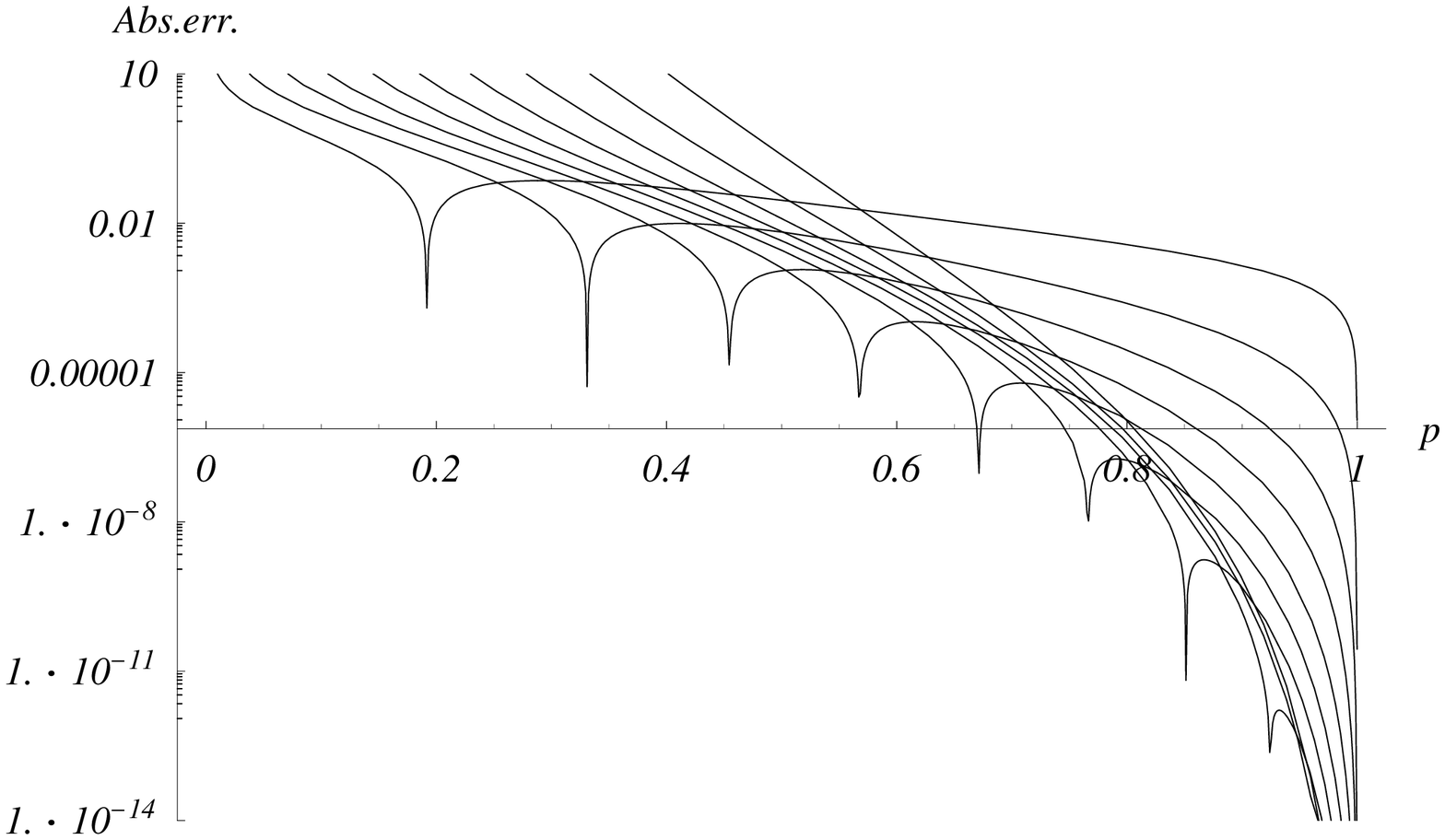}
\includegraphics[width=16cm]{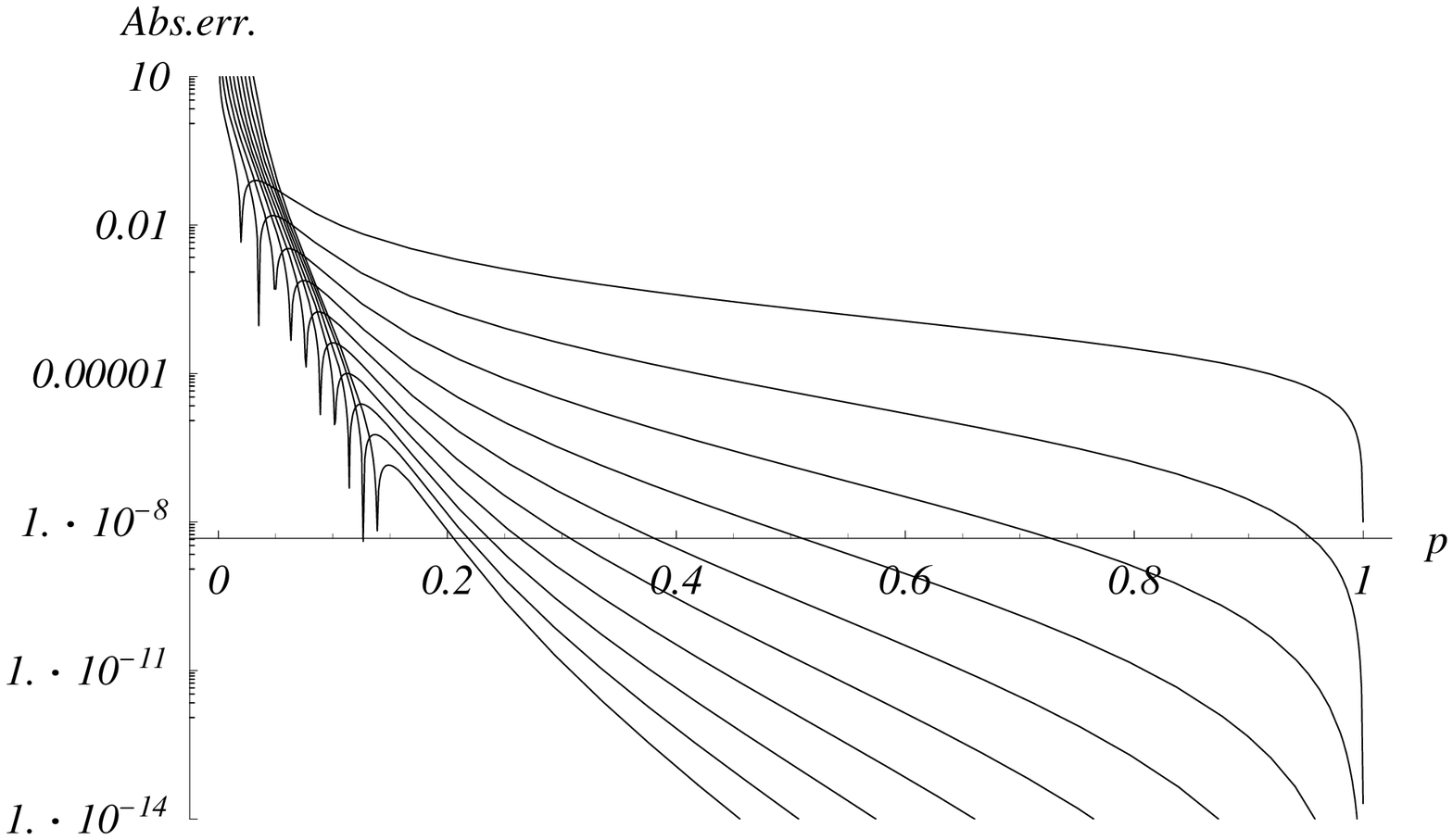}
\caption{
Absolute error as a function of $p$ of Rempala's expansion of
the first inverse moment of a positive binomial variate for $N=10$ (upper graph) and $N=100$ (lower graph),
with 1 term (upper curve), and up to 10 terms (lowest curve).
\label{fig:rempala10}
}
\end{figure}

A recently obtained expansion is Rempala's \cite{Rempala:03}, which is an improvement on an expansion
by Marciniak and Weso{\l}owski \cite{Marciniak:99}:
\be
(1-q^N) \E[1/K|K>0] = (Np)^{-1} \sum_{i=0}^{M-1} \frac{(q/p)^i}{{N-1\choose i}}.
\ee
For not too small $p$ this series seems to converge much faster than Stephan's expansion.
However, as it is an asymptotic expansion, it actually diverges.
For fixed $p$ there is an optimal number of terms $M$,
and for larger $M$ the error increases very rapidly. For example, taking $N=100$,
the first 100 terms of the series give extremely accurate results for $p>0.54$
but are completely useless for $p<0.51$.
The absolute error of this expansion is illustrated in Figure \ref{fig:rempala10}.

\begin{figure}
\includegraphics[width=16cm]{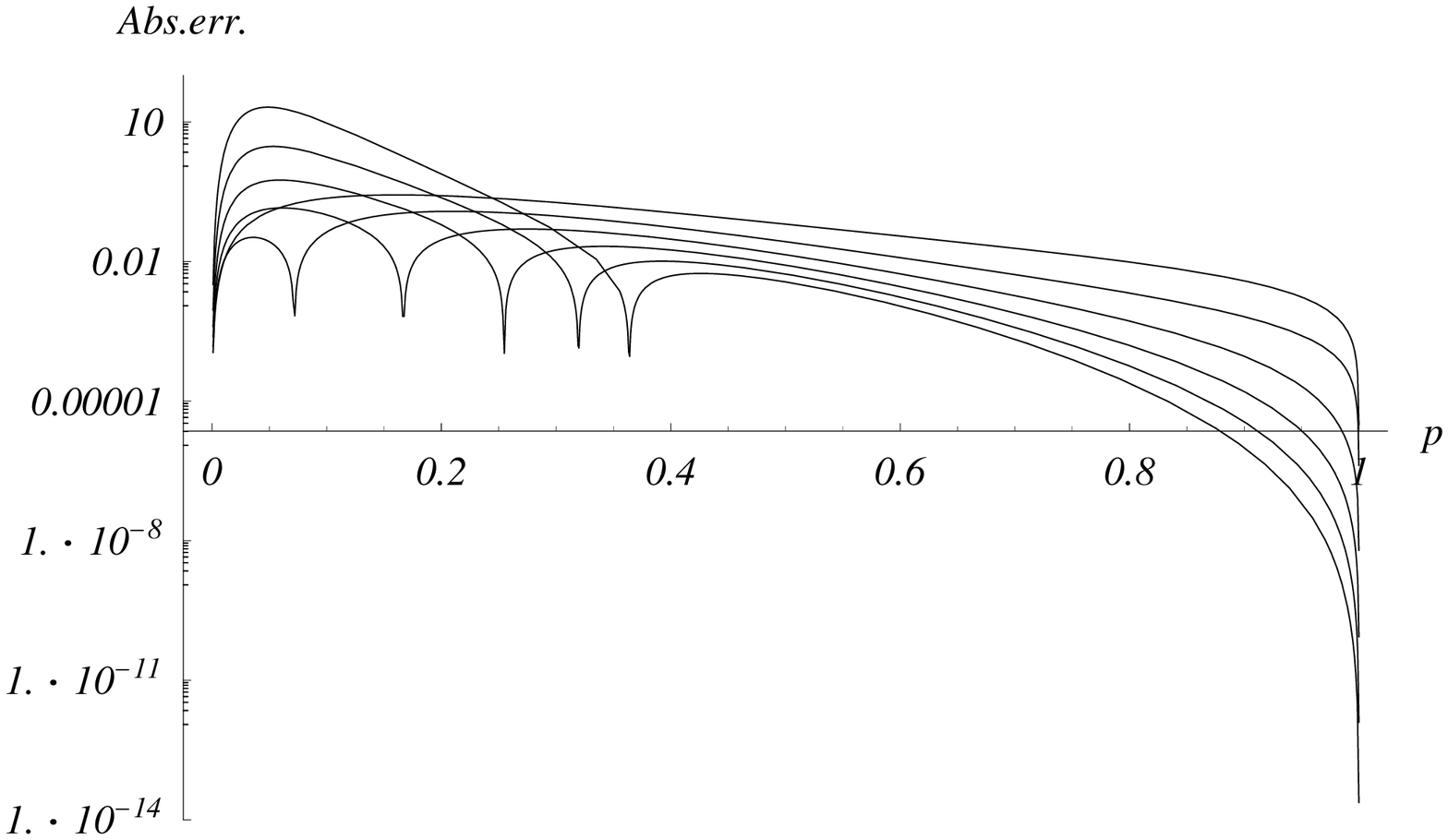}
\includegraphics[width=16cm]{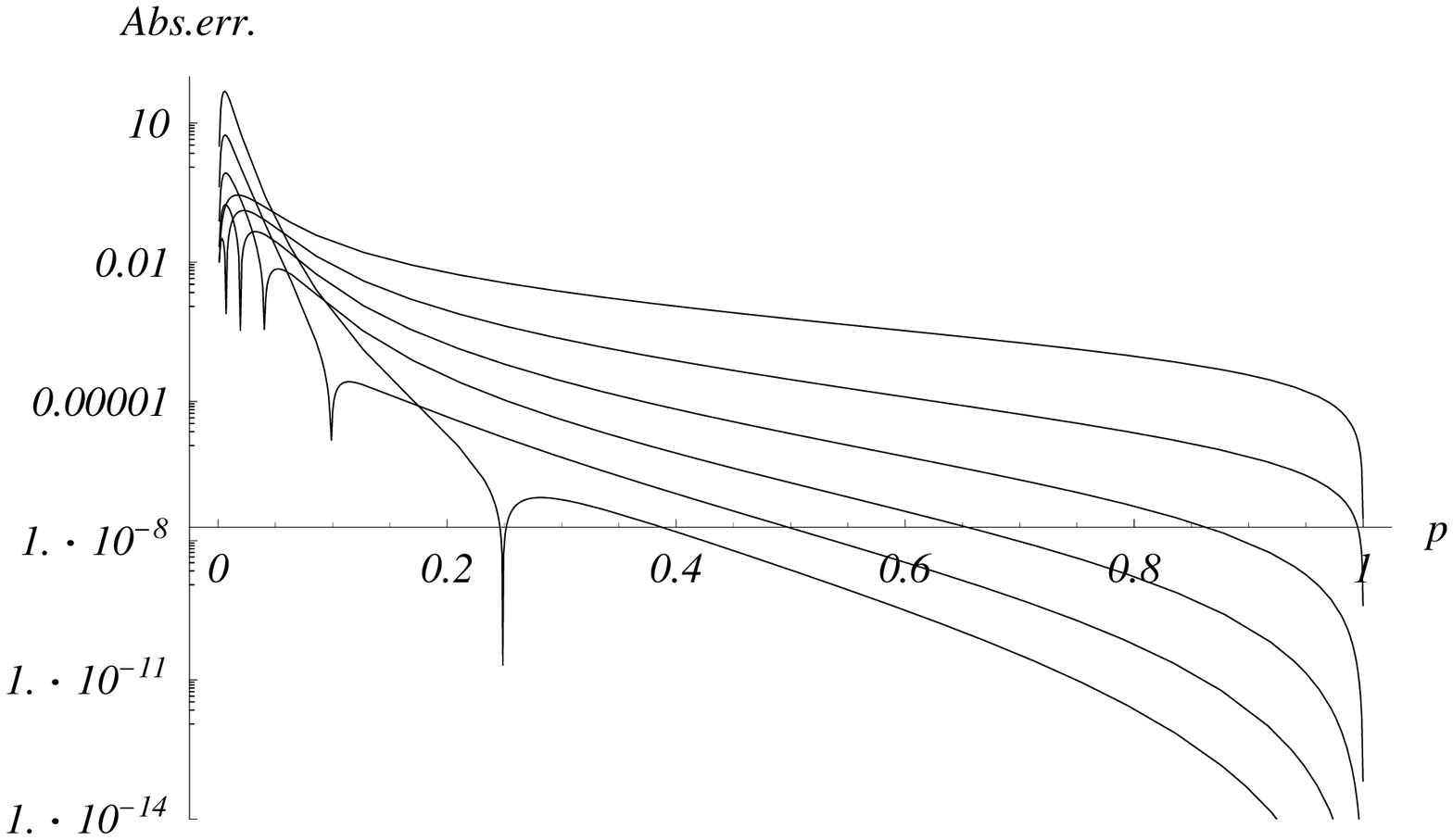}
\caption{
Absolute error as a function of $p$ of Znidaric's expansion of
the first inverse moment of a positive binomial variate for $N=10$ (upper graph) and $N=100$ (lower graph),
with 1 term (upper curve), 3,5,7,9 and 11 terms (lowest curve).
Expansions with an even number of terms are left out for clarity.
\label{fig:znid10}
}
\end{figure}

Znidaric \cite{Znidaric:05} gives an expansion formula for all inverse moments.
When specialised to the case of the first inverse moment it reads:
\be
(1-q^N) \E[1/K|K>0] = \frac{Np}{(Np+q)^2}\,\sum_{i=0}^{M-1}(-1)^i(i+1)\frac{\mu_i(N-1)}{(Np+q)^i},
\ee
where $\mu_i(N-1)$ is the $i$-th central moment of $\mbox{Bin}(N-1,p)$.
Just like Rempala's expansion, this is an asymptotic series and suffers from the same divergence problems;
in addition it gives much less accurate results, as seen from Figure \ref{fig:znid10}.

\section{Poisson expansion of a probability density\label{sec:barbour}}
It is well-known that the binomial distribution $\mbox{Bin}(N,p)$,
with probability distribution $f(k) = {N\choose k}p^k(1-p)^{n-k}$,
$k=0,\ldots,N$, tends to the Poisson distribution
with mean value $\mu=Np$, when $Np$ is kept fixed and $N$ tends to infinity.
In \cite{Barbour:87,Barbour:89} Barbour and coworkers presented a series expansion
for a general distribution of a non-negative discrete random variate
in terms of the Poisson distribution, providing a quantitative and more generally useful
version of this last statement. They gave this expansion the name of Poisson-Charlier expansion.

The first term of the Poisson expansion of a distribution $f$ is, of course, the Poisson distribution itself,
with mean $\mu$ given by the mean of $f$:
\be
\mu = \E[k] = \sum_{k=0}^\infty f(k) k.
\ee
The higher-order terms of the Poisson expansion consist of backward differences ($\nabla$) 
of $\pi_\mu(k)$, with coefficients based on the factorial cumulants $\kappa^{(k)}$ of $f$.
Formally, the Poisson expansion is based on the following identity:
\be
f = \exp\left(\sum_{k=2}^\infty \frac{\kappa^{(k)}}{k!}(-\Nabla)^k\right) \pi_\mu.
\label{eq:poissonidentity}
\ee
A simple proof of this identity is given in Appendix \ref{app:proof}.

In practice, the right-hand side of (\ref{eq:poissonidentity}) 
has to be replaced by a finite series, with a finite number of terms in
powers of $\nabla$, and with suitable error bounds estimating the truncation error.
There are many ways to do this.
The most obvious way is to replace the right-hand side by a Taylor series in $\nabla$.
An $m$-th order series can be obtained by expanding the function
\be
t\mapsto \exp\left(\sum_{k=2}^\infty \frac{\kappa^{(k)}}{k!}(-t\Nabla)^k\right)
\label{eq:seriesalt}
\ee
in powers of $t$ up to $t^{m-1}$ around $t=0$, and then putting $t=1$.
This yields a polynomial in $\nabla$ of degree $m-1$.
We will generally define the \textit{order} of the expansion as 1 plus the highest degree of $t$,
and its \textit{degree} as the highest degree of $\nabla$.
While conceptually simple, we are not aware of known error estimates for this kind of expansion.

In contrast, Barbour's expansion amounts to an $m$-th order series of the function
\be
t\mapsto \exp\left(\frac{1}{t}\sum_{k=2}^\infty \frac{\kappa^{(k)}}{k!}(-t\Nabla)^k\right),
\label{eq:operator}
\ee
which has an extra factor $1/t$ in the argument of the exponential function.
Taylor expansion in $t$ yields a polynomial of degree $2(m-1)$ in $\nabla$, 
hence the degree of the $m$-th order expansion is $2(m-1)$. 
When comparing series (\ref{eq:seriesalt}) and (\ref{eq:operator}) of the same degree one finds that
the latter contains a subset of the terms of the former.
For example, the third order series (\ref{eq:operator}) is given by
\be
1+\frac{\kappa^{(2)}}{2}\,\Nabla^2
 -\frac{\kappa^{(3)}}{6}\,\Nabla^3
 +\frac{(\kappa^{(2)})^2}{8}\,\Nabla^4,
\label{eq:poissonexp}
\ee
while the fifth order series (\ref{eq:seriesalt}), which is also of degree 4, has an additional term:
\be
1+\frac{\kappa^{(2)}}{2}\,\Nabla^2
-\frac{\kappa^{(3)}}{6}\,\Nabla^3
+\left(\frac{(\kappa^{(2)})^2}{8}+\frac{\kappa^{(4)}}{24}\right)\,\Nabla^4.
\ee
Surprisingly, however, while this observation makes the series (\ref{eq:seriesalt}) look more powerful
than (\ref{eq:operator}) (more terms for the same degree), the latter is much better suited to the purposes 
of this paper because it actually converges much faster and, moreover, explicit error bounds for it are known.

The complete expression of (\ref{eq:operator}), 
including all high-order terms, is quite complicated and will not be given here
(see \cite{Barbour:87}, eq.~(2.7)).
Applying this operator polynomial to $\pi_\mu(k)$ yields the $m$-term approximation of $f(k)$.
General bounds on the absolute error are given in \cite{Barbour:87}. 
These bounds were obtained in a highly non-trivial way (not to say magical way), 
using the so-called Stein's method \cite{Barbour:05}.
We will not give the general bounds here but will only mention them (in Section \ref{sec:main}) 
for the special case of the binomial distribution.
 
\medskip

\noindent\textbf{Example 1.} 
The Poisson expansion of the Poisson distribution trivially reduces to the first term only, as it should,
because all factorial cumulants of the Poisson distribution are 0, except $\kappa^{(1)}$ which is
equal to the mean value $\mu$.

\medskip

\noindent\textbf{Example 2.}
A more enlightening example is the Poisson expansion of the binomial distribution
$K\sim\mbox{Bin}(N,p)$.
Its mean is $\mu=Np$ and $\E[(1+x)^K] = (1+px)^N$, so that
$$
\log\E[(1+x)^K] = N\log(1+px) = N\,\sum_{j=1}^\infty \frac{(-1)^{j+1}}{j}(px)^j.
$$
Hence, the $j$-th factorial cumulant (for $j>0$) is
\be
\kappa^{(j)} = -N(j-1)!(-p)^j.\label{eq:bincum}
\ee
Thus, the third order Poisson expansion of the PDF of $K$ is
\be
f(k)\approx
\pi_{\mu}(k)
-\frac{\mu^2}{2N}\,\Nabla^2 \pi_{\mu}(k)
-\frac{\mu^3}{3N^2}\,\Nabla^3 \pi_{\mu}(k)
+\frac{\mu^4}{8N^3}\,\Nabla^4 \pi_{\mu}(k)
\ee
with $\mu=Np$.
\section{Inverse moments of the Poisson distribution\label{sec:expectation}}
Using the Poisson expansion, calculating the $r$-th inverse moment of a distribution is reduced
to calculating the expectations $\E[1/(Q+a)^r|Q>0]$ for $a=0,1,\ldots$ where $Q$ is a positive
Poisson variate with mean value $\mu$.
The expectations for $a=0$ are essentially the inverse moments of the positive Poisson distribution.
The other expectations are the inverse moments about $-a$ and can be derived from the moments about the
origin using a simple recurrence relation.

To simplify notations we introduce the symbol $\E^+$:
for any function $g$, and for $K$ a positive random variate with probability distribution $f$:
\beas
\E^+[g(K)] &:=& \sum_{k=1}^\infty f(k) g(k) \\
&=& (1-f(0))\E[g(K)|K>0] \\
&=& \E[g(K)] - f(0)g(0).
\eeas
Thus, in particular,
\be
\E^+[1/(Q+a)^r] = \E[1/(Q+a)^r]-e^{-\mu}/a^r.\label{eq:exppos}
\ee

Let's consider first the inverse moments about the origin,
$\E^+[1/Q^r]$.
For $r=1$, we have
\bea
\E^+[1/Q] &=& e^{-\mu}\,\sum_{i=1}^\infty \frac{\mu^i}{i\, i!} \label{eq:eiseries1} \\
&=& e^{-\mu}\,\int_0^1 \md t\,\, \frac{e^{\mu t}-1}{t}. \label{eq:eiintegral}
\eea
The integral is essentially equal to the exponential integral $\mbox{Ei}$ \cite{Abramowitz:72}. With
$\gamma\approx 0.5772\ldots$ the Euler-Mascheroni constant,
\bea
\int_0^1 \md t \,\,\frac{e^{\mu t}-1}{t}
&=& \mbox{Ei}(\mu)-\log\mu-\gamma=:\mbox{Er}(\mu). \label{eq:erdef}
\eea
The exponential integral function is well-studied, and implementations are incorporated in many
numerical and algebraic software packages.

The higher inverse moments about the origin can be expressed in terms of hypergeometric functions
\be
\E^+[1/Q^r] = \mu e^{-\mu}\,\, \mbox{}_{r+1}F_{r+1}\left({1\,1\,\ldots1 \atop 2\,2\,\ldots2}\,\Bigg|\,\mu\right),
\ee
which is actually just a standardised restatement of the definition.
Again, these functions can be accurately and efficiently calculated using standard software packages.
In absentia, a second option is to resort to explicit series expansions;
some of these are given in Appendix \ref{app:numerical}.


\bigskip

Now we move on to the (non-central) inverse moments about $-a$, $\E[1/(Q+a)^r]$, for $a>0$.
Let's first consider the case $r=1$.
For $a=1$,
\be
\E[1/(Q+1)] = \frac{1-e^{-\mu}}{\mu},\label{eq:rec0}
\ee
and for $a>1$ we can use the recurrence \cite{Chao:72}
\be
\E[1/(Q+a)] = \frac{1}{\mu}\left(1-(a-1)\E[1/(Q+a-1)]\right). \label{eq:rec1}
\ee
This recurrence can be generalised to higher $r$ as follows:
\begin{lemma}
For $Q$ a positive Poisson variate with mean $\mu$, and for integers $a>0$ and $r>1$:
if $a=1$,
\be
\E[1/(Q+1)^r] = \frac{1}{\mu}\E^+[1/Q^{r-1}], \label{eq:rec2}
\ee
and if $a>1$,
\be
\E[1/(Q+a)^r] = \frac{1}{\mu}\left(\E[1/(Q+a-1)^{r-1}]-(a-1)\E[1/(Q+a-1)^r]\right). \label{eq:rec3}
\ee
\end{lemma}
In \cite{Jones:04} a similar recurrence formula is used.

\noindent\textit{Proof.}
First note
\beas
\E[1/(Q+a)^r] &=& e^{-\mu}\sum_{k=0}^\infty \frac{\mu^k}{k!}\,\frac{1}{(k+a)^r}\\
&=& e^{-\mu}\sum_{k=0}^\infty \frac{\mu^k}{(k+1)!}\,\frac{k+1}{(k+a)^r}\\
&=& \frac{1}{\mu}e^{-\mu}\sum_{k=0}^\infty \frac{\mu^{k+1}}{(k+1)!}\,\frac{k+1}{(k+a)^r}\\
&=& \frac{1}{\mu}e^{-\mu}\sum_{k=1}^\infty \frac{\mu^{k}}{k!}\,\frac{k}{(k+a-1)^r} \\
&=& \frac{1}{\mu}\E^+[Q/(Q+a-1)^r].
\eeas
For $a=1$, this directly gives (\ref{eq:rec2}).

For $a>1$, the $k=0$ term is zero and can be added to the sum:
\beas
\E[1/(Q+a)^r] &=& \frac{1}{\mu}e^{-\mu}\sum_{k=0}^\infty \frac{\mu^{k}}{k!}\,\frac{k}{(k+a-1)^r} \\
&=& \frac{1}{\mu}\E[Q/(Q+a-1)^r].
\eeas
Writing $Q = (Q+a-1) - (a-1)$ then yields
\beas
\E[1/(Q+a)^r] &=& \frac{1}{\mu}\E[Q/(Q+a-1)^r] \\
&=& \frac{1}{\mu}\left(\E[1/(Q+a-1)^{r-1}]-(a-1)\E[1/(Q+a-1)^r]\right).
\eeas
\qed

\bigskip

By iterating the recurrence relations (\ref{eq:rec0}), (\ref{eq:rec1}), (\ref{eq:rec2}) and (\ref{eq:rec3}),
one can find an expression of the $r$-th inverse moments of the positive Poisson distribution
in terms of central and non-central Stirling numbers.
\begin{proposition}
Let $Q$ be a positive Poisson variate with mean value $\mu$.
For integers $a\ge1$ and $r\ge1$:
\be
\E[1/(Q+a)^r] = \frac{1}{\mu^a}\left(
S_a^{(r)}\,(1-e^{-\mu}) + \sum_{k=1}^{r-1} S_a^{(r-k)}\,\mu_{-k}
+\sum_{k=1}^{a-1} S_{a-k,k}^{(r)}\,\mu^k
\right).\label{eq:ekar}
\ee
\end{proposition}
\textit{Proof.}
Let's denote the right hand side of (\ref{eq:ekar}) by $A(a,r)$, and the $r$-th inverse moment of the positive Poisson distribution $\E^+[1/Q^r]$ by $\mu_{-r}$.
We consider first the case $r=a=1$.
Since $A(1,1) = (1-e^{-\mu})/\mu$ it coincides with $\E[1/(Q+1)]$, by (\ref{eq:rec0}).
Next, for $r=1$ and $a>1$, we need to check the recurrence (\ref{eq:rec1}), i.e.\ whether
$$
A(a,1) = (1-(a-1) A(a-1,1))/\mu.
$$
From the generating function (\ref{eq:genshift}), we have
$S_{a-k,k}^{(1)} = (-1)^{a-k-1} (a-1)!/k!$. Thus,
\beas
A(a,1) &=& \frac{(-1)^{a-1}(a-1)!}{\mu^a}\left(1-e^{-\mu}+\sum_{k=1}^{a-1}(-\mu)^k/k!\right) \\
&=& \frac{(-1)^{a-1}(a-1)!}{\mu^a}\left(-\sum_{k=a}^{\infty}(-\mu)^k/k!\right).
\eeas
It is straightforward to check the recurrence from the latter expression.

For $r>1$ and $a=1$ we find
$A(1,r)=(\sum_{k=1}^{r-1} S_{1}^{(r-k)}\mu_{-k}+S_1^{(r)}(1-e^{-\mu}))/\mu$.
As $S_1^{(r)} = 1$ if and only if $r=1$, this simplifies to $\mu_{-(r-1)}/\mu$, as required for (\ref{eq:rec2}).

Finally, to check the remaining recurrence (\ref{eq:rec3}) for $a,r>1$, we first note that the $k=a-1$ sum term
in the last term of (\ref{eq:ekar}) vanishes, because $S_{1,a}^{(r)}=0$ for $r>1$. Hence the upper summation
limit can be replaced by $a-2$. Then, from recurrences (\ref{eq:recstir}) and (\ref{eq:recshift}), we see
that the coefficients appearing in the three summation terms obey one and the same recurrence:
$S_{a}^{(r-k)} = S_{a-1}^{(r-k-1)}-(a-1)S_{a-1}^{(r-k)}$,
$S_{a}^{(r)} = S_{a-1}^{(r-1)}-(a-1)S_{a-1}^{(r)}$, and
$S_{a-k,k}^{(r)} = S_{a-1-k,k}^{(r-1)}-(a-1)S_{a-1-k,k}^{(r)}$.
It is then an easy matter to verify that $A(a,r)$ indeed satisfies the final recurrence
$A(a,r)=A(a-1,r-1)-(a-1)A(a-1,r)$.

Since $A(a,r)$ satisfies the same recurrences and boundary conditions of $\E[1/(Q+a)^r]$, the two must coincide,
proving equality in (\ref{eq:ekar}).
\qed

\bigskip

Inserting the explicit formulas (\ref{eq:stirling1}) and (\ref{eq:stirlingnc1}) 
immediately gives, for the case $r=1$:
\be
\E[1/(Q+a)] = \frac{(a-1)!(-1)^{a-1}}{\mu^a}\left(1-e^{-\mu}+\sum_{j=1}^{a-1}(-\mu)^j/j!\right).
\label{eq:e1ka}
\ee
\section{Poisson Expansion of Inverse Moments\label{sec:main}}
Based on the results of the two previous Sections, we are now in the position to formulate our main result.
While our own interest lies with the binomial distribution, our result is generally valid for any
positive discrete random variate.
\begin{theorem}
Let $K$ be a positive discrete random variate with probability distribution $f$, having mean value $\mu$ and
factorial cumulants $\kappa^{(i)}$.

Let $\mu_{-r}$ be the $r$-th inverse moment $\mu_{-r} = \E^+[1/Q^r]$
of a positive Poisson variate $Q$ with the same mean value $\mu$ as $K$,
and let $\mu_{-r,a}$ be the corresponding shifted inverse moments $\mu_{-r,a} = \E[1/(Q+a)^r]$, for $a>0$,
obtainable e.g.\ from (\ref{eq:ekar}).
Define the sequence $q_{-r}$ with $q_{-r}(a)=\mu_{-r,a}$, for $a\ge0$.

Let $P_m$ be the degree $2(m-1)$ polynomial obtained from the $m$-th order
Taylor approximation of (\ref{eq:operator}) as indicated in Section \ref{sec:barbour}.

Then the $m$-th order Poisson approximation of the $r$-th inverse moment of $K$ is given by
\be
\E^+[1/K^r] \approx (P_m(-\Delta) q_{-r})(0).\label{eq:theorem1}
\ee
\end{theorem}
\textit{Proof.}
Since the Poisson expansion of a distribution is expressed in terms of $\Nabla^j \pi_\mu$, we first
calculate the following sums:
\beas
\nu_{-r,j} &:=& \sum_{k=1}^\infty \frac{1}{k^r}\,\Nabla^j \pi_\mu(k) \\
&=& \sum_{k=1}^\infty \frac{1}{k^r}\,\sum_{a=0}^j {j \choose a}(-1)^a \pi_\mu(k-a) \\
&=& \sum_{k=1}^\infty \frac{1}{k^r} \pi_\mu(k)
  + \sum_{a=1}^j {j \choose a}(-1)^a \sum_{k=a}^\infty \frac{\pi_\mu(k-a)}{k^r}.
\eeas
Now note
\beas
\sum_{k=a}^\infty \frac{\pi_\mu(k-a)}{k^r} &=& \sum_{k=0}^\infty \frac{\pi_\mu(k)}{(k+a)^r} \\
&=& \E[1/(Q+a)^r].
\eeas
Thus
\beas
\nu_{-r,j} &=& \E^+[1/Q^r] + \sum_{a=1}^j {j \choose a}(-1)^a \E[1/(Q+a)^r] \\
&=& \sum_{a=0}^j {j \choose a}(-1)^a \mu_{-r,a}.
\eeas
We can express the last equation in terms of the sequence $q_{-r}$ and 
the $j$-th forward difference operator:
$$
\nu_{-r,j} = ((-\Delta)^j q_{-r})(0).
$$
Combining this with the $m$-term Poisson expansion of the distribution of $K$,
$$
f \approx P_m(\Nabla) \pi_\mu,
$$
gives the final result
$$
\E^+[1/K^r] \approx (P_m(-\Delta) q_{-r})(0).
$$
\qed

\bigskip

\begin{figure}[ht]
\includegraphics[width=16cm]{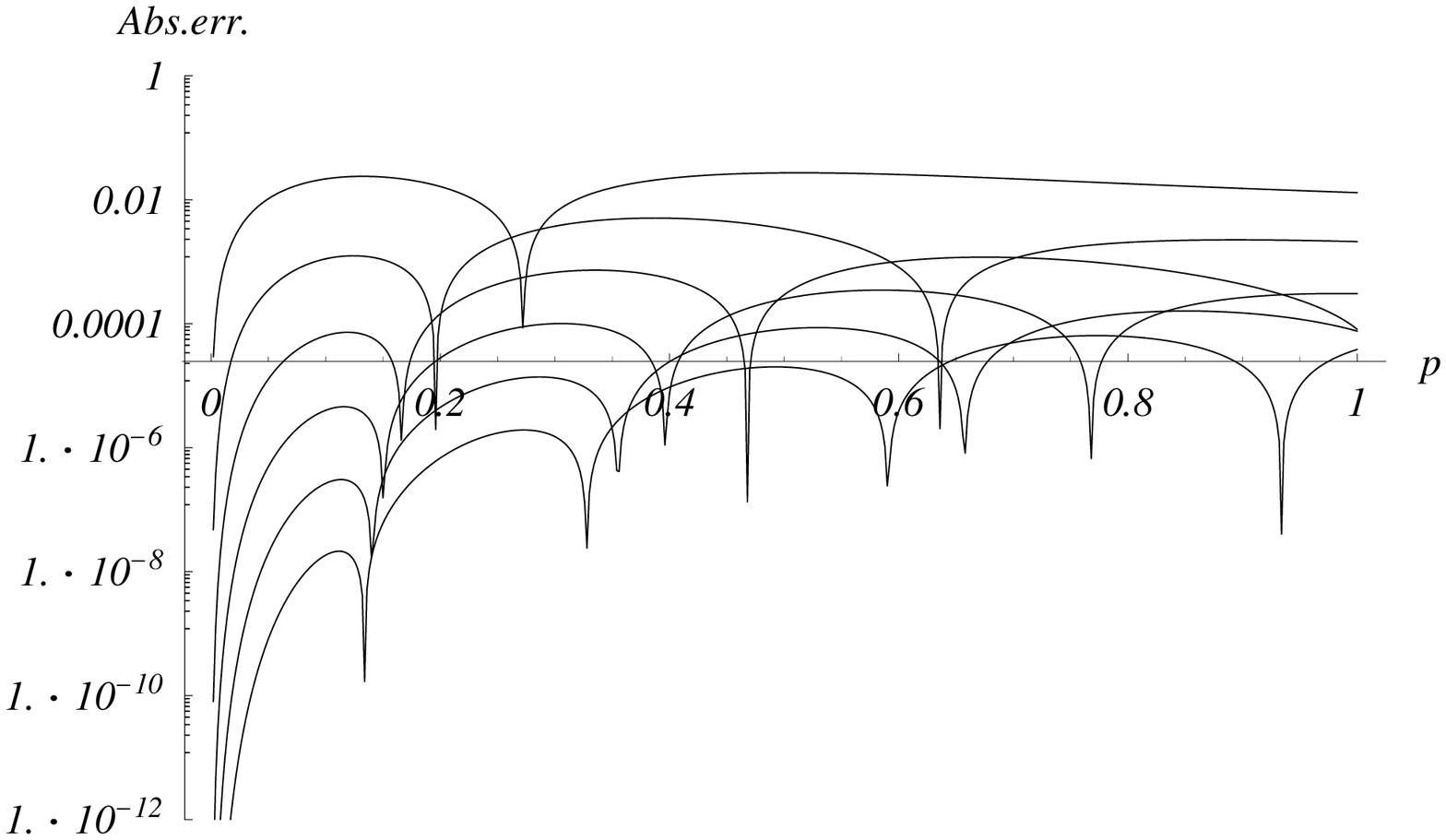}
\includegraphics[width=16cm]{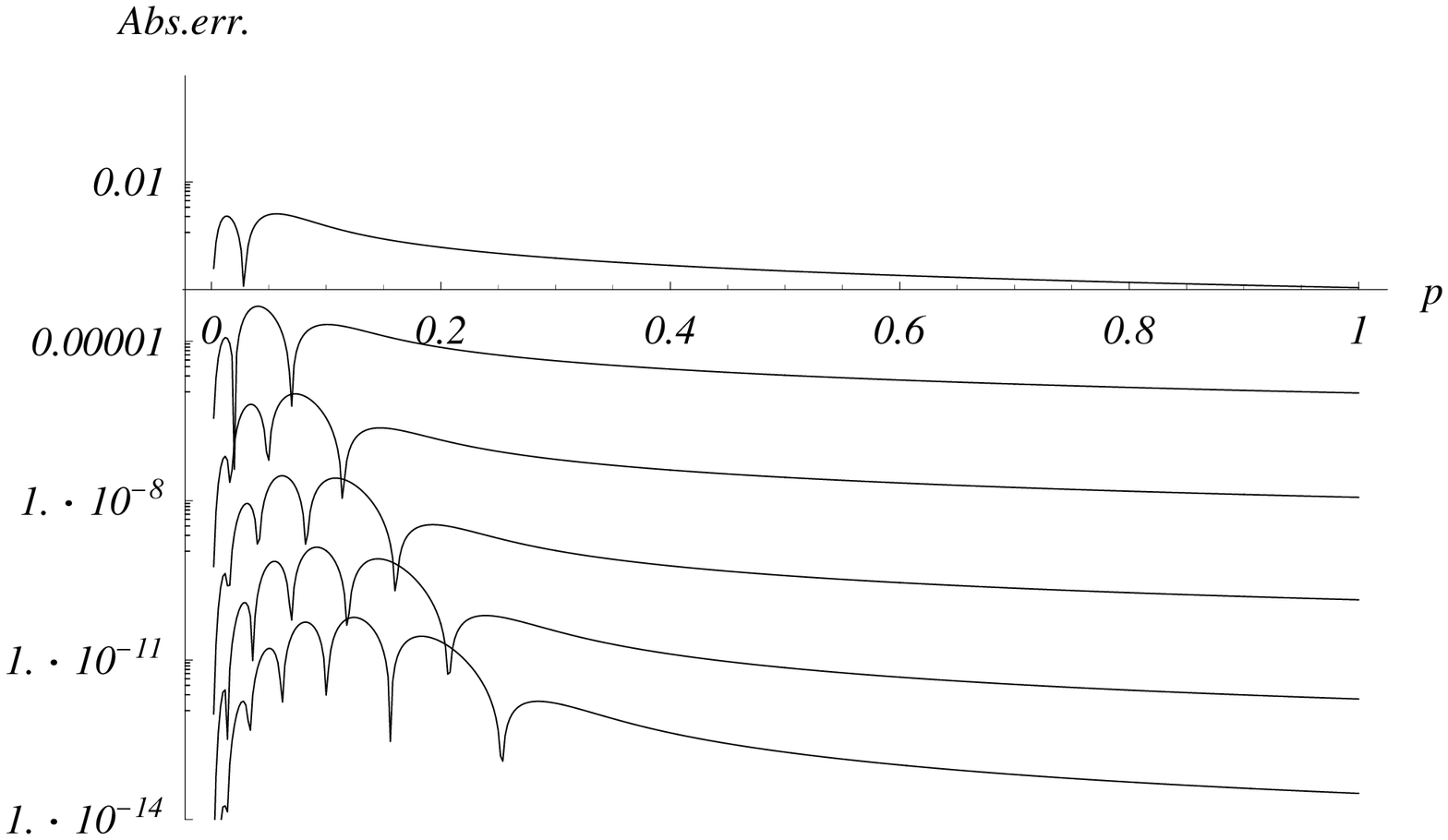}
\caption{
Absolute error as a function of $p$ of the Poisson expansion of
the first inverse moment of a positive binomial variate for $N=10$ and $N=100$,
with 1 term (upper curve), and up to 6 terms (lowest curve).
\label{fig:poisson10abs}
}
\end{figure}

\begin{figure}[ht]
\includegraphics[width=16cm]{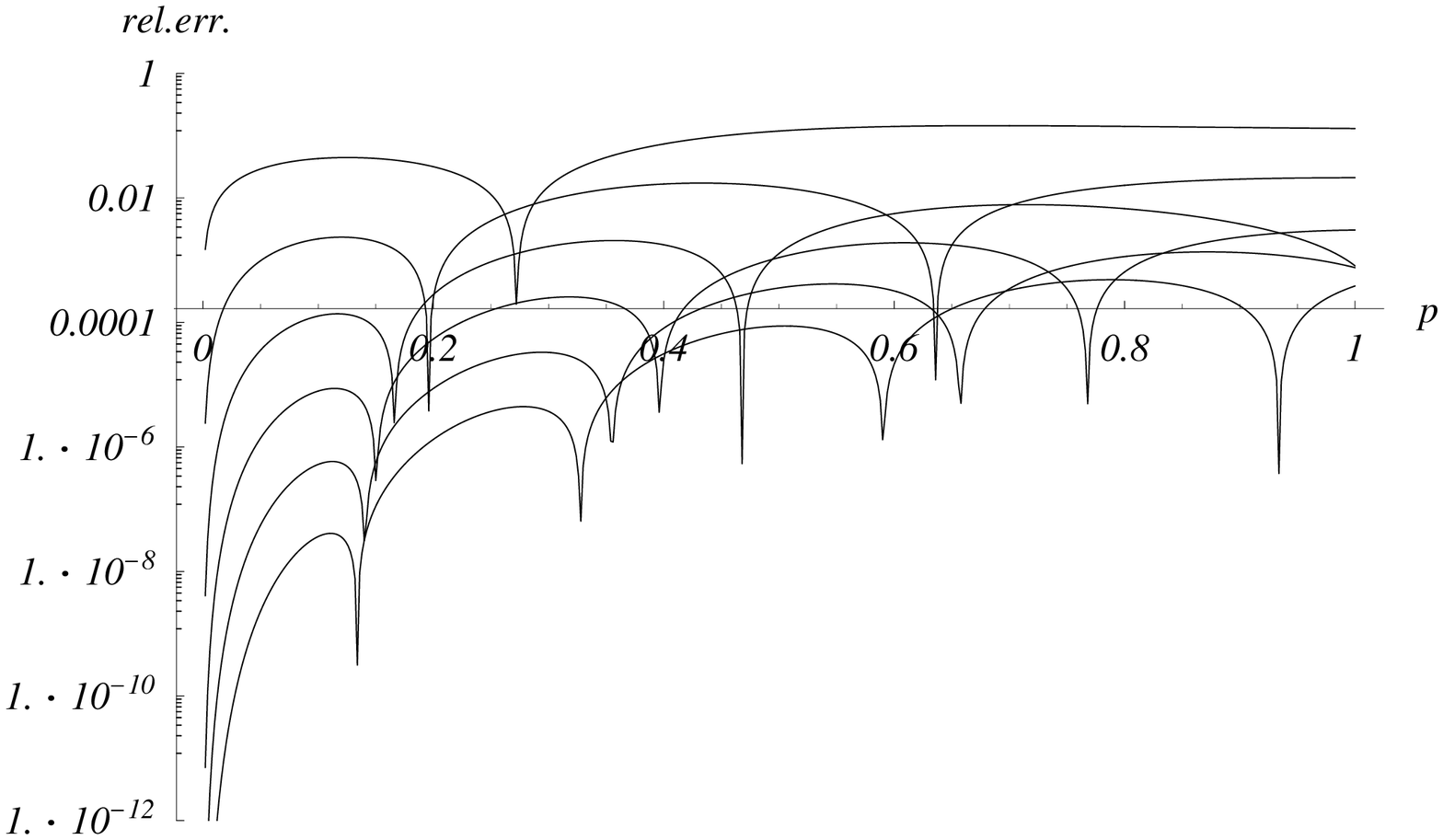}
\includegraphics[width=16cm]{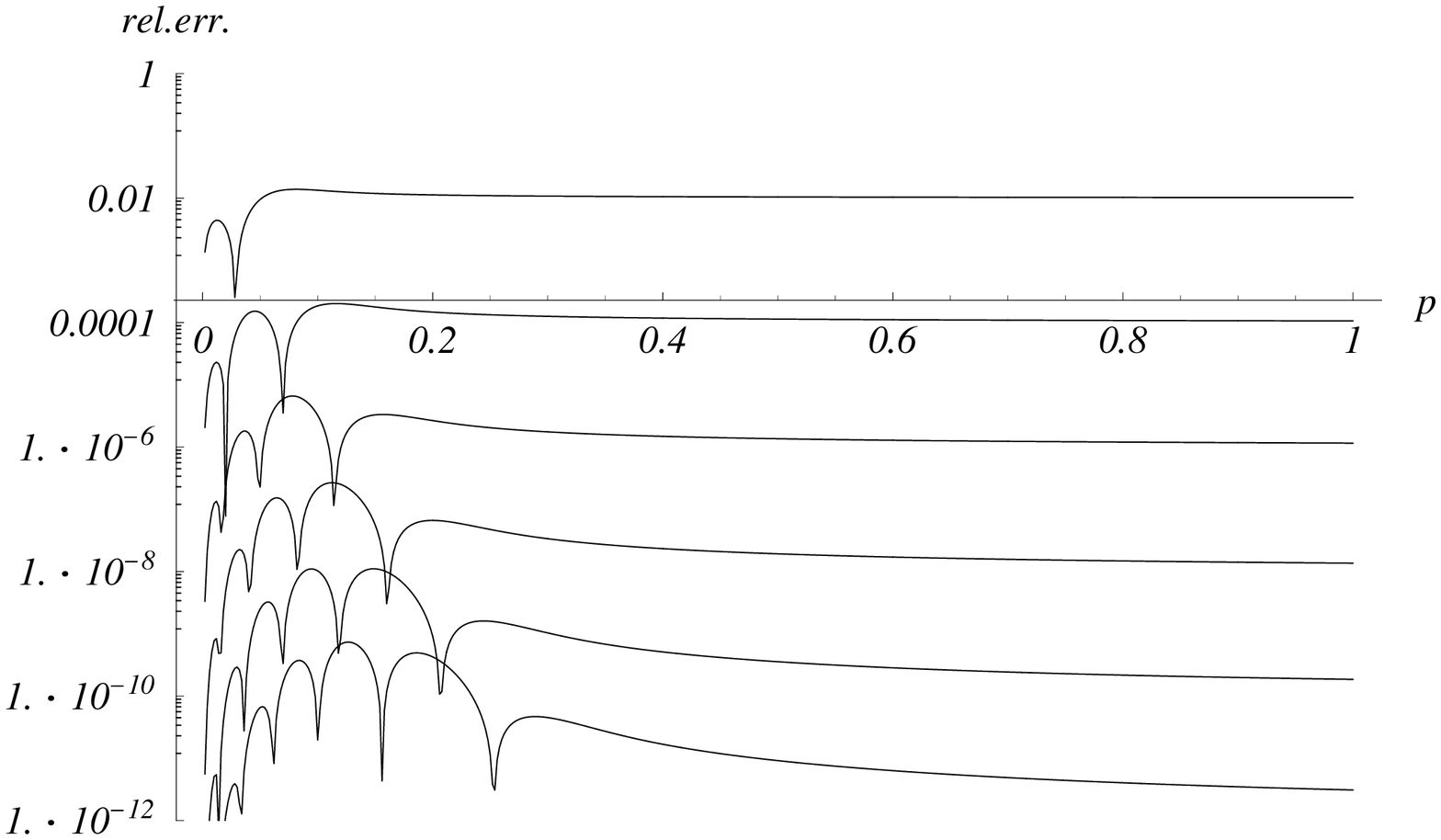}
\caption{
Same as Figure \ref{fig:poisson10abs} but showing the relative error.
\label{fig:poisson10}
}
\end{figure}

When applied to the positive binomial distribution, this Theorem yields the following expansion (substituting
formula (\ref{eq:bincum}) for the factorial cumulants):
\begin{corollary}
For $K$ a positive binomial variate $K\sim \mbox{\rm Bin}(N,p)$, the $m$-th Poisson approximation
of the $r$-th inverse moment is given by
\be
\E^+[1/K^r] \approx 
q_{-r}(0) + \sum_{k=1}^{m-1} \frac{1}{N^k}\sum_{j=1}^k \frac{(-1)^j}{j!} \,\alpha_{k-j,j}\,\mu^{j+k}
((-\Delta)^{j+k}q_{-r})(0).\label{eq:corollary1}
\ee
Here, the coefficients $\alpha_{l,j}$ obey the recurrence
\be
\alpha_{l,j+1} = \sum_{k=0}^l \frac{\alpha_{k,j}}{l-k+2},
\ee
with boundary conditions $\alpha_{0,0}=1$ and $\alpha_{l,0}=0$ for $l>0$.
\end{corollary}
Note that when fixing $\mu$, this expansion is a series expansion in $1/N$.

\bigskip

\noindent\textit{Proof.}
With $\kappa^{(j)} = -N(j-1)!(-p)^j$, (\ref{eq:operator}) becomes
\beas
t &\mapsto& \exp\left( -\frac{N}{t}\,\sum_{k=2}^\infty \frac{1}{k} (t p \Nabla)^k\right) \\
&=& \sum_{j=0}^\infty \frac{1}{j!} \,
\left( -\frac{N}{t} \,\sum_{k=2}^\infty \frac{1}{k} (t p \Nabla)^k \right)^j \\
&=& \sum_{j=0}^\infty \frac{1}{j!} \, (-N)^j \sum_{l=0}^\infty \alpha_{l,j}t^{j+l}(p\Nabla)^{2j+l}.
\eeas
Here we have defined the series expansion
$$
\left(\sum_{k=2}^\infty \frac{x^k}{k}\right)^j = \sum_{l=0}^\infty \alpha_{l,j}x^{2j+l}.
$$
It is easily checked from this definition
that the coefficients $\alpha_{l,j}$ obey the recurrence stated in the corollary.
The $m$-term approximation polynomial $P_m$ is now obtained by imposing the constraint
$j+l\le m-1$ and setting $t=1$, giving 
$$
P_m(x) = \sum_{j=0}^{m-1} \frac{1}{j!} \, (-N)^j
\sum_{l=0}^{m-1-j} \alpha_{l,j}\,(px)^{2j+l}.
$$
Substituting $p=\mu/N$ and collecting terms in identical powers of $N$ gives
\beas
P_m(x) &=& 1 + \sum_{j=1}^{m-1} \frac{(-1)^j}{j!} 
\sum_{l=0}^{m-1-j} \alpha_{l,j}\,N^{-j-l} (\mu x)^{2j+l} \\
&=& 1+\sum_{k=1}^{m-1} \frac{1}{N^k}\sum_{j=1}^k \frac{(-1)^j}{j!} \,\alpha_{k-j,j}\,(\mu x)^{j+k}.
\eeas
To obtain the last line we have set $k=j+l$ and reorganised the double summation.
Combining this formula with (\ref{eq:theorem1}) then gives the formula of the corollary.
\qed

\bigskip

Table \ref{tab:alpha} gives the first values of $\alpha_{l,j}$.
One sees that $\alpha_{l,1}=1/(l+2)$ and $\alpha_{l,2} = 2(H_{l+2}-1)/(l+4)$, 
where $H_n$ is the $n$-th harmonic number.
We are not aware of any closed form expression for $j>2$.
\begin{table}[h]
\begin{tabular}{l|cccccccc}
$l\backslash j$ & 0 &1& 2& 3& 4& 5& 6& 7 \\
\hline
0&1&$1/2$&$1/4$&$1/8$&$1/16$&$1/32$&$1/64$&$1/128$ \\
1&0&$1/3$&$1/3$&$1/4$&$1/6$&$5/48$&$1/16$&\\
2&0&$1/4$&$13/36$&$17/48$&$7/24$&$125/576$&&\\
3&0&$1/5$&$11/30$&$59/135$&$229/540$&&&\\
4&0&$1/6$&$29/80$&$241/480$&&&&\\
5&0&$1/7$&$223/630$&&&&&\\
6&0&$1/8$&&&&&&\\
7&0&&&&&&&
\end{tabular}
\caption{Values of the coefficients $\alpha_{l,j}$ used in Corollary 1, for $j+l\le7$.
\label{tab:alpha}}
\end{table}

For the special case $r=1$, we can present an even more explicit formula.
\begin{corollary}
For $K$ a positive binomial variate $K\sim \mbox{\rm Bin}(N,p)$, the $m$-th Poisson approximation
of its first inverse moment is given by
\be
\E^+[1/K] \approx 
y_0 + \sum_{k=1}^{m-1} \frac{1}{N^k} \left(
\sum_{j=1}^k \frac{(-1)^j}{j!}\,\alpha_{k-j,j}\,y_{j+k}
\right),
\ee
where
\be
y_n = \mu^n e^{-\mu}\mbox{\rm Er}(\mu) + \sum_{l=1}^n (l-1)!\left(e^{-\mu}{n\choose l}-1\right)\mu^{n-l},
\ee
and $\mu=Np$.
\end{corollary}
\textit{Proof.}
Equations (\ref{eq:eiintegral}), (\ref{eq:erdef}) and (\ref{eq:e1ka}) 
yield explicit formulas for the sequence $q_{-1}$:
\beas
\mu_{-1,0} &=& e^{-\mu} \mbox{Er}(\mu) \\
\mu_{-1,a} &=& \frac{(a-1)!(-1)^{a-1}}{\mu^a}\left(1-e^{-\mu}+\sum_{j=1}^{a-1}(-\mu)^j/j!\right).
\eeas
Then this gives
\beas
((-\Delta)^{n}q_{-1})(0)
&=& e^{-\mu}\mbox{Er}(\mu)
-\sum_{a=1}^n {n\choose a}\frac{(a-1)!}{\mu^a}\left(-e^{-\mu}+\sum_{j=0}^{a-1}(-\mu)^j/j!\right)\\
&=& e^{-\mu}\mbox{Er}(\mu)
+e^{-\mu}\sum_{a=1}^n {n\choose a}\frac{(a-1)!}{\mu^a}
-\sum_{a=1}^n {n\choose a}\frac{(a-1)!}{\mu^a}\sum_{j=0}^{a-1}(-\mu)^j/j!
\eeas
The third term simplifies, upon setting $l=a-j$ and rearranging the double sum:
\beas
\sum_{a=1}^n {n\choose a}\frac{(a-1)!}{\mu^a}\sum_{j=0}^{a-1}(-\mu)^j/j!
&=& \sum_{a=1}^n\sum_{j=0}^{a-1} {n\choose a}\frac{(a-1)!}{j!}(-1)^j \mu^{j-a} \\
&=& \sum_{l=1}^n \mu^{-l} \sum_{j=0}^{n-l} {n\choose j+l}\frac{(j+l-1)!}{j!}(-1)^j \\
&=& \sum_{l=1}^n (l-1)!\mu^{-l}.
\eeas
That yields
$$
((-\Delta)^{n}q_{-1})(0) = 
e^{-\mu}\mbox{Er}(\mu)
+\sum_{l=1}^n \left(e^{-\mu}{n\choose l}-1\right) \frac{(l-1)!}{\mu^l}.
$$
Substituting this expression for $((-\Delta)^{n}q_{-1})(0)$ in (\ref{eq:corollary1}) 
of Corollary 1 gives the desired result.
\qed

\bigskip

To illustrate the behaviour of the expansion, we depict the absolute and relative error of the approximation
of $\E^+[1/K]$ in Figures \ref{fig:poisson10abs} and \ref{fig:poisson10}, respectively, for $N=10$ and $N=100$.
It is clear from these figures that in contrast to previous expansions, the error is
bounded uniformly over the complete interval $[0,1]$.
The graphs have been produced using a Mathematica program, listed in Appendix \ref{app:mathematica}.
We have also calculated the absolute error for the alternative series expansion (\ref{eq:seriesalt}). 
It turned out that this expansion converged much more slowly than a Barbour expansion of the same degree 
(let alone one of the same order),
even though the latter contains fewer terms than the former. For want of a better explanation, we attribute this
phenomenon to the magic of Stein's method.

The results depicted in Figure \ref{fig:poisson10abs} can be compared to some 
explicit error bounds in \cite{Barbour:87}.
Corollary 2.4 in \cite{Barbour:87} gives an upper bound to the absolute error when approximating expectations
of a sum $W$ of $N$ Bernoulli (`0-1') variates $X_i$. With $p_i=P[X_i=1]$, the absolute error
of the $m$-th order approximation to the expectation $\E[h(W)]$ is bounded as
$$
|\eta_m| \le 2^{2m-1}\frac{1-e^{-\mu}}{\mu}\sum_i p_i^{m+1}||h||
$$
where $\mu=\E[W] = \sum_i p_i$.
A binomial variate $K\sim\mbox{Bin}(N,p)$ is just a special case of this, obtained by taking all $p_i$ equal.
For the inverse moments, $h$ is given by $h=(0,1,2^{-r},3^{-r},\ldots)$, so that $||h||=1$.
This gives the following bound:
\bea
|\eta_m| &\le& 2^{2m-1}\frac{1-e^{-\mu}}{\mu}N p^{m+1} \nonumber \\
&=& 2^{2m-1}(1-e^{-\mu})p^m.
\eea
Obviously, this bound is only useful for $p<1/4$, and it matches the actual convergence only for $p\ll 1$.
Nevertheless, this bound partially proves our claim that the Poisson approximation to the inverse moments 
converges. For bigger values of $p$, we currently have to rely on the numerical calculations reported in
Figure \ref{fig:poisson10abs}.
%
\appendix
\section{Numerical calculation of the inverse moments of a positive Poisson variate}
\label{app:numerical}
\begin{figure}
\includegraphics[width=12cm]{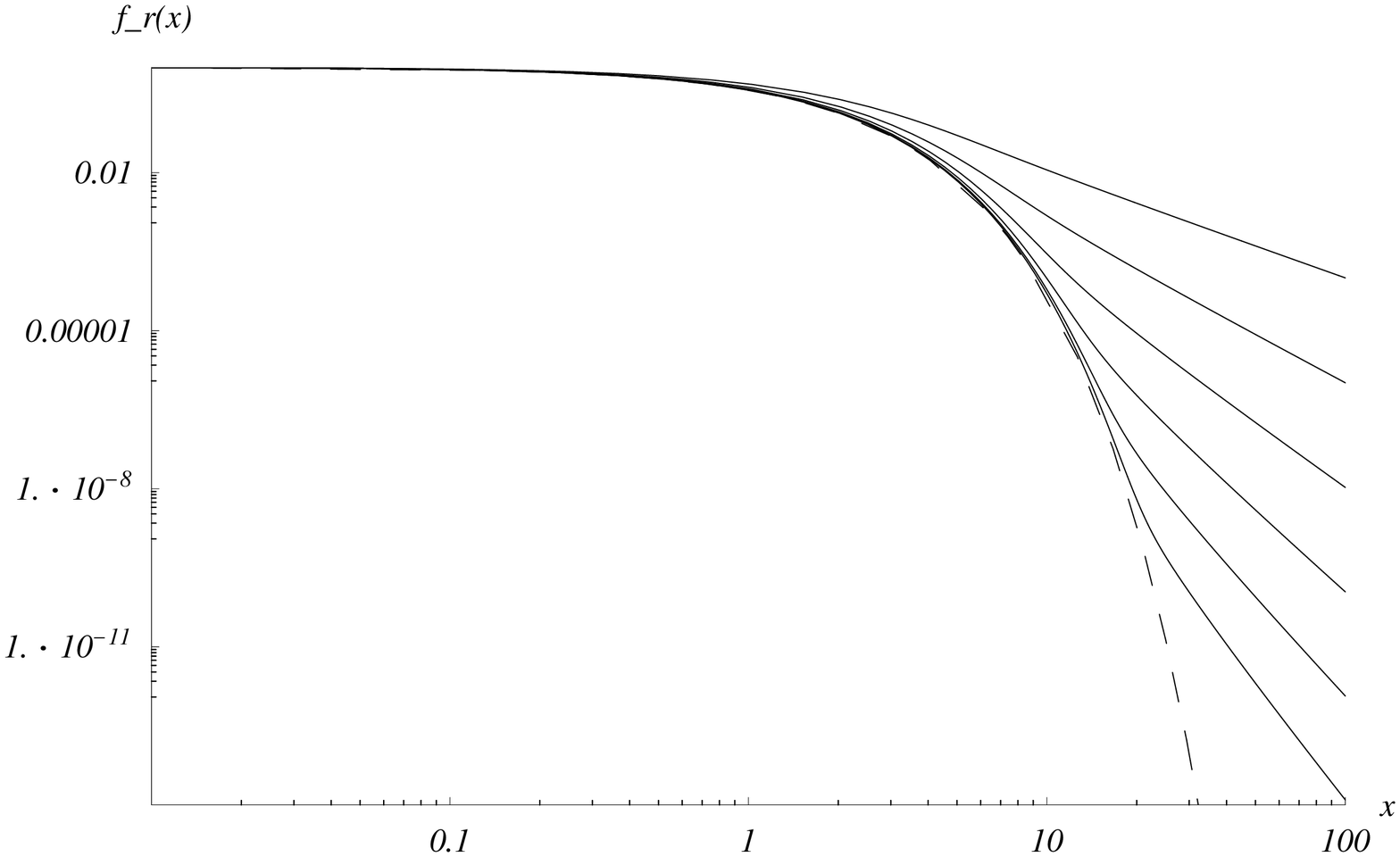}
\caption{
Figure of $f_r(\mu)/\mu=\frac{e^{-\mu}}{\mu}\,\sum_{k=1}^\infty \frac{\mu^k}{k!\,k^r}$,
for values of $r$ from 1 (upper curve) to 6 (lower curve).
The dashed curve represents the function $e^{-\mu}$.
\label{fig:poissonmoments}
}
\end{figure}

In this Appendix we present a numerical method for calculating the inverse moments
$f_r(\mu) = \E^+[1/Q^r]$ of a positive Poisson variate with mean value $\mu$.
Series expansions are given in \cite{Gupta:79,Jones:04,Stancu:68,Tiku:64,Znidaric:05} but as these are asymptotic
series they are not universally applicable.
Plotting $f_r(\mu)/\mu$ reveals two different regimes in the range of $\mu$.
For small $\mu$ this function is seen to behave as $e^{-\mu}$ for all $r$ (see Figure \ref{fig:poissonmoments}).
This is also clear from the definition of $f_r$, as the first terms of its defining series are
$$
\frac{f_r(\mu)}{\mu} = e^{-\mu}(1+\frac{\mu}{2.2!} +\frac{\mu^2}{3.3!}+\ldots).
$$
This suggests that for small values of $\mu$, the truncated series should give a good approximation:
\be
f_r(\mu) \approx e^{-\mu}\,\sum_{k=1}^{M_1}\frac{\mu^k}{k!\,k}.
\label{eq:pos1}
\ee
For larger $\mu$, the Figure suggests a $1/\mu^r$ behaviour, with a moderately sharp cross-over region.
For these larger values we can use the asymptotic series of \cite{Jones:04}:
\be
f_r(\mu) \approx \sum_{i=0}^{M_2-1}\frac{|s_{r+i}^{(r)}|}{\mu^{r+i}},
\label{eq:pos2}
\ee
where $s_j^{(k)}$ are the Stirling numbers of the first kind.

By appropriately choosing the cross-over point $\mu_*$ at which to switch
from (\ref{eq:pos1}) to (\ref{eq:pos2})
and the number of terms $M_1$ and $M_2$ in the two series,
one can tune the maximal relative error of the approximation while keeping the computational effort at bay.
Tables \ref{tab:app} and \ref{tab:app10} show the values for the cross-over point $\mu_*$,
and $M_1$ and $M_2$, as function of $r$,
needed to obtain an approximation with relative error ($|$1-approximation/exact value$|$) below
$10^{-5}$ and $10^{-10}$, respectively.
\begin{table}[h]
\begin{tabular}{rllllll}
\hline
$r$     &  1 & 2 & 3 & 4 & 5 & 6 \\ \hline
$\mu_*$ & 13.671 & 17.061 & 20.544 & 24.775 & 28.966 & 32.969 \\
$M_1$   & 31     & 35     & 39     & 44     & 49     & 53 \\
$M_2$   & 10     & 15     & 20     & 26     & 32     & 38 \\ \hline
\end{tabular}
\caption{Values for the cross-over point $\mu_*$, and optimal number of terms $M_1$ and $M_2$, as function of $r$,
to obtain an approximation of $f_r(\mu)$ with relative error below $10^{-5}$.\label{tab:app}}
\end{table}
\begin{table}[h]
\begin{tabular}{rllllll}
\hline
$r$     &  1 & 2 & 3 & 4 & 5 & 6 \\ \hline
$\mu_*$ & 25.734 & 29.206 & 33.998 & 37.903 & 42.573 & 47.068 \\
$M_1$   & 63     & 67     & 74     & 79     & 85 & 90 \\
$M_2$   & 20     & 26     & 33     & 39     & 46 & 53 \\ \hline
\end{tabular}
\caption{Same as Table \ref{tab:app} but for a relative error below $10^{-10}$.\label{tab:app10}}
\end{table}
In general, the choice of values for $M_1$ and $M_2$ involves a trade-off
between the two series. In the present case, however, the choice of $M_2$ is determined because the second series
is an asymptotic expansion. From a certain number of terms onwards, the benefit in including additional terms
becomes marginal and ultimately the series diverges. We have chosen the value of $M_2$ that minimises
the value $\mu_*$ below which the relative error becomes larger than the set minimum, so that the number of terms
$M_1$ of the first series, covering the remaining interval, can be made as small as possible.
\section{Proof of the Poisson Expansion Formula}
\label{app:proof}
In this Appendix we give a simple proof of the identity (\ref{eq:poissonidentity}) underlying the Poisson expansion.
Let $f$ be a semi-infinite sequence
$$
f=(f(0),f(1),\ldots,f(k),\ldots).
$$
As everywhere in this paper, we set
$f(k)=0$ for $k<0$.
Let $f$ have mean value $\mu$ and let its $k$-th factorial cumulant be $\kappa^{(k)}$,
with generating function
$$
\log\left(\sum_{k=0}^\infty f(k)(1+x)^k\right) = \sum_{j=0}^\infty \frac{\kappa^{(j)}}{j!}\,x^j.
$$
Recall that for any PDF $f$, $\kappa^{(0)}=0$ and $\kappa^{(1)}=\mu$.

Define the operator $S$
$$
S:=\exp\left(\sum_{k=2}^\infty \frac{\kappa^{(k)}}{k!}\,(-\nabla)^k\right),
$$
having matrix representation
$$
S\mapsto \bm{S} = \exp\left(\sum_{k=2}^\infty \frac{\kappa^{(k)}}{k!}\,(-\bm{\nabla})^k\right).
$$
Let $g$ denote the sequence $g=S \pi_\mu= \bm{S}.\pi_\mu$, where $\pi_\mu$ is the sequence
$$
(\pi_\mu(0),\pi_\mu(1),\ldots,\pi_\mu(k),\ldots)
$$
of the PDF of the Poisson distribution with mean value $\mu$.
We thus need to prove that $g=f$.

\medskip

The generating function of factorial cumulants can be written as
$$
\log\left(\sum_{k=0}^\infty f(k)(1+x)^k\right) = \log f^T\xi(x),
$$
where $\xi(x)$ is the semi-infinite vector
$$
\xi(x) = (1,1+x,(1+x)^2,\ldots,(1+x)^k,\ldots),
$$
for $x\in\C$.
This vector is the eigenvector of $\bm{\Delta}$ corresponding to eigenvalue $x$.
Applying $\Delta$ to the sequence $\xi(x)$ indeed yields $\Delta\xi(x) = (x,x(1+x),x(1+x)^2,\ldots) = x\xi(x)$.

We will now calculate the factorial cumulant generating function of $g$.
The inner product $g^T\xi(x)$ is given by
\beas
g^T\xi(x) &=& \pi_\mu^T. \bm{S}^T. \xi(x) \\
&=& \pi_\mu^T.\exp\left(\sum_{k=2}^\infty \frac{\kappa^{(k)}}{k!}\,(-\bm{\nabla})^k\right)^T.\xi(x) \\
&=& \pi_\mu^T.\exp\left(\sum_{k=2}^\infty \frac{\kappa^{(k)}}{k!}\,\bm{\Delta}^k\right).\xi(x).
\eeas
Since $\xi(x)$ is an eigenvector of $\bm{\Delta}$ with eigenvalue $x$, this immediately gives
\beas
g^T\xi(x) 
&=& \exp\left(\sum_{k=2}^\infty \frac{\kappa^{(k)}}{k!}\,x^k\right)\,\,\pi_\mu^T.\xi(x).
\eeas
The logarithm of the last factor is the factorial cumulant generating function of the Poisson distribution,
which is known to be $\mu x$.
Thus we get that the factorial cumulant generating function of $g$ is
\beas
\log g^T\xi(x) &=& \sum_{k=2}^\infty \frac{\kappa^{(k)}}{k!}\,x^k + \mu x \\
&=& \sum_{k=0}^\infty \frac{\kappa^{(k)}}{k!}\,x^k.
\eeas
Since the right-hand side is identical to the factorial cumulant generating function of $f$,
we have proven that $f=g$.
\qed
\section{A Mathematica program for the first inverse moment of a positive binomial variate}
\label{app:mathematica}
Here we reproduce the short Mathematica program that we have used to prepare Figure \ref{fig:poisson10abs}.
For ease of implementation, the inverse moments of the Poisson distribution are calculated directly using 
Mathematica's summation capabilities, rather than via any recurrences or explicit formulas 
like the one of Proposition 1.

\begin{verbatim}
(* First inverse moment of a-shifted Poisson: *)
invmom[mu_,a_] := Sum[Exp[-mu]mu^k/k!/(k+a),{k,If[a==0,1,0],Infinity}]

(* First inverse moment of l-th forward difference of Poisson: *)
invdif[mu_,l_] := Sum[Binomial[l,j](-1)^j invmom[mu,j],{j,0,l}]

(* Factorial cumulants of binomial: *)
kappa[n_, p_, k_] = -n(k - 1)!(-p)^k;

(* m term Poisson expansion :*)
expansion[n_, mu_, x_, m_] :=
    Collect[Normal[
          Series[Exp[Sum[kappa[n, mu/n, k](-x t)^k/k!, {k, 2, m}]/t], 
          {t, 0, m - 1}]] /. t -> 1, x];
              
(* The m-th order Poisson approximation (m=1,...,6) to 
   the first inverse moment of Bin(n, p) :*)
appr[n_, p_] =
    Table[(expansion[n, n p, x, m] /. x^k_ -> invdif[n p, k]) - 1 +
        invdif[n p, 0], {m, 1, 6}];
        
(* Exact expression: *)
exact[n_, p_] = Sum[Binomial[n, k]p^k(1 - p)^(n - k)/k, {k, 1, n}];

(* Absolute error: *)
abserr[p_, n_, m_] := Abs[appr[n, p][[m]]-exact[n, p]]

(* Relative error: *)
relerr[p_, n_, m_] := Abs[1 - appr[n, p][[m]]/exact[n, p]]

(* Produces the graph of Fig. 4 for n=10, m = 1 to 6: *)
<< Graphics`Graphics`
grlist = Table[
      LinearLogListPlot[
        Table[{p, N[relerr[SetPrecision[p, 30], 10, m], 30]} /.
            p -> k/500, {k, 1, 500}], PlotRange -> All,
        PlotJoined -> True], {m, 1, 6}];
Show[grlist, PlotRange -> {-12, 0}, DefaultFont -> {"Times-Italic", 16},
    AxesLabel -> {"p", "rel.err."}];
\end{verbatim}

\end{document}